\documentclass[12pt]{article}%
\usepackage{sw20bams}
\usepackage{amsmath}
\usepackage{amsfonts}
\usepackage{amssymb}
\usepackage{graphicx}%
\providecommand{\U}[1]{\protect\rule{.1in}{.1in}}
\begin{document}

\title{Median Area for Broken Sticks}
\author{Steven Finch}
\date{April 25, 2018}
\maketitle

\begin{abstract}
Breaking a line segment $L$ in two places at random, the three pieces can be
configured as a triangle $T$ with probability $1/4$. \ We determine both the
PDF\ and CDF\ for $\operatorname*{area}(T)$ in terms of elliptic integrals.
\ In particular, if $L$ has length $1$, then the median area $0.031458...$ can
be calculated to arbitrary precision. \ We also mention the analog involving
cyclic quadrilaterals -- with corresponding probability $1/2$ -- and ask some
unanswered questions.

\end{abstract}

For simplicity's sake, we start with a stick of length $2$ (not $1$). \ A
triangle, formed by two independent uniform breaks in the stick, has sides
$a+b+c=2$ satisfying%
\[%
\begin{array}
[c]{ccccc}%
0<a<b+c, &  & 0<b<a+c, &  & 0<c<a+b
\end{array}
\]
hence%
\[%
\begin{array}
[c]{ccccc}%
0<a<1, &  & 0<b<1, &  & 1<a+b<2.
\end{array}
\]
The joint density function for $(a,b)$ is thus $2$ (constant) over the shaded
triangular region in Figure 1; the marginal density for $a$ is $2x$ if
$0<x<1$; the cross-correlation between $a$ and $b$ is $-1/2$. \
%TCIMACRO{\FRAME{ftbpFU}{2.9983in}{3.0234in}{0pt}{\Qcb{Triangular support for
%bivariate side density}}{}{triangle.eps}%
%{\special{ language "Scientific Word";  type "GRAPHIC";
%maintain-aspect-ratio TRUE;  display "USEDEF";  valid_file "F";
%width 2.9983in;  height 3.0234in;  depth 0pt;  original-width 4.9978in;
%original-height 5.0401in;  cropleft "0";  croptop "1";  cropright "1";
%cropbottom "0";  filename 'triangle.eps';file-properties "XNPEU";}} }%
%BeginExpansion
\begin{figure}[ptb]%
\centering
\includegraphics[
height=3.0234in,
width=2.9983in
]%
{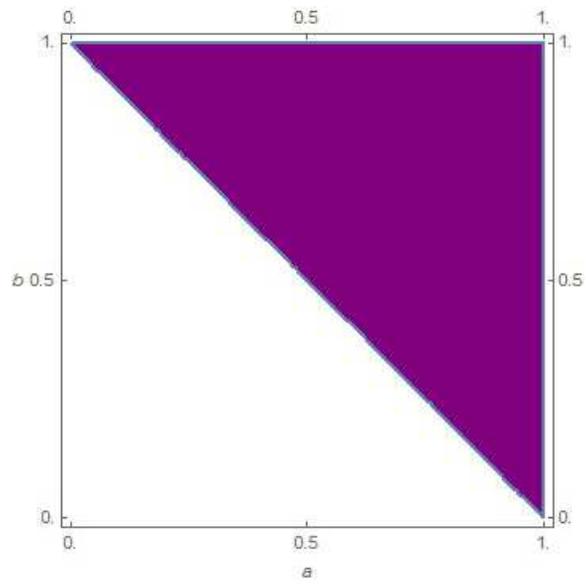}%
\caption{Triangular support for bivariate side density}%
\end{figure}
%EndExpansion
By Heron's formula, the mean area of the triangle is \cite{Dobbs}%
\[
\mathbb{E}\left(  \operatorname*{area}\right)  =%
%TCIMACRO{\dint \limits_{0}^{1}}%
%BeginExpansion
{\displaystyle\int\limits_{0}^{1}}
%EndExpansion%
%TCIMACRO{\dint \limits_{1-x}^{1}}%
%BeginExpansion
{\displaystyle\int\limits_{1-x}^{1}}
%EndExpansion
2\sqrt{(1-x)(1-y)(x+y-1)}\,dy\,dx=\frac{4\pi}{105}%
\]
and the mean square area is%
\[
\mathbb{E}\left(  \operatorname*{area}\nolimits^{2}\right)  =%
%TCIMACRO{\dint \limits_{0}^{1}}%
%BeginExpansion
{\displaystyle\int\limits_{0}^{1}}
%EndExpansion%
%TCIMACRO{\dint \limits_{1-x}^{1}}%
%BeginExpansion
{\displaystyle\int\limits_{1-x}^{1}}
%EndExpansion
2(1-x)(1-y)(x+y-1)\,dy\,dx=\frac{4^{2}}{960}=\frac{1}{60}.
\]

In contrast, a cyclic quadrilateral \cite{Johnson, CxtrGrtz}, formed by three
independent uniform breaks in the stick, has sides $a+b+c+d=2$ in this order
satisfying%
\[%
\begin{array}
[c]{ccccccc}%
0<a<b+c+d, &  & 0<b<a+c+d, &  & 0<c<a+b+d, &  & 0<d<a+b+c
\end{array}
\]
hence%
\[%
\begin{array}
[c]{ccccccc}%
0<a<1, &  & 0<b<1, &  & 0<c<1, &  & 1<a+b+c<2.
\end{array}
\]
The joint density for $(a,b,c)$ is thus $3/2$ (constant) over the shaded
hexahedral region in Figure 2; note the additional complexity of two missing
corners, not just one. \
%TCIMACRO{\FRAME{ftbpFU}{5.6317in}{3.0242in}{0pt}{\Qcb{Hexahedral support for
%trivariate side density (two views)}}{}{hexahedra.eps}%
%{\special{ language "Scientific Word";  type "GRAPHIC";
%maintain-aspect-ratio TRUE;  display "USEDEF";  valid_file "F";
%width 5.6317in;  height 3.0242in;  depth 0pt;  original-width 11.2253in;
%original-height 6.0027in;  cropleft "0";  croptop "1";  cropright "1";
%cropbottom "0";  filename 'hexahedra.eps';file-properties "XNPEU";}} }%
%BeginExpansion
\begin{figure}[ptb]%
\centering
\includegraphics[
height=3.0242in,
width=5.6317in
]%
{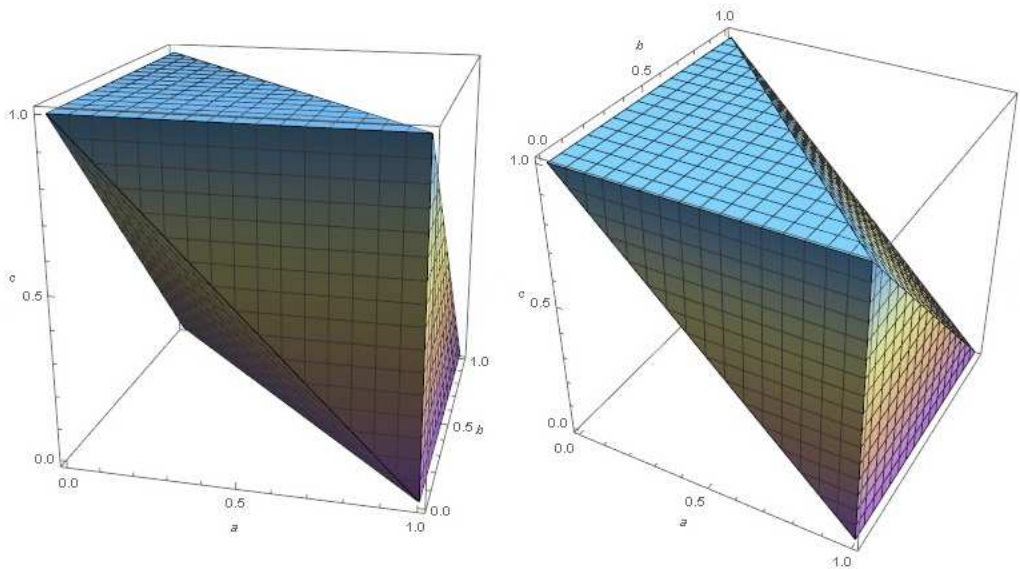}%
\caption{Hexahedral support for trivariate side density (two views)}%
\end{figure}
%EndExpansion
The marginal density for $a$ is $\frac{3}{4}\left(  1+2x-2x^{2}\right)  $ for
$0<x<1$; the cross-correlation between $a$ and $b$ is $-1/3$. \ By
Brahmagupta's formula, the mean area of the cyclic quadrilateral is
\cite{PACM1}
\begin{align*}
\mathbb{E}\left(  \operatorname*{area}\right)   &  =%
%TCIMACRO{\dint \limits_{0}^{1}}%
%BeginExpansion
{\displaystyle\int\limits_{0}^{1}}
%EndExpansion
\;%
%TCIMACRO{\dint \limits_{0}^{1-x}}%
%BeginExpansion
{\displaystyle\int\limits_{0}^{1-x}}
%EndExpansion
\,%
%TCIMACRO{\dint \limits_{1-x-y}^{1}}%
%BeginExpansion
{\displaystyle\int\limits_{1-x-y}^{1}}
%EndExpansion
\frac{3}{2}\sqrt{(1-x)(1-y)(1-z)(x+y+z-1)}\,dz\,dy\,dx\\
&  +%
%TCIMACRO{\dint \limits_{0}^{1}}%
%BeginExpansion
{\displaystyle\int\limits_{0}^{1}}
%EndExpansion
\;%
%TCIMACRO{\dint \limits_{1-x}^{1}}%
%BeginExpansion
{\displaystyle\int\limits_{1-x}^{1}}
%EndExpansion
\;%
%TCIMACRO{\dint \limits_{0}^{2-x-y}}%
%BeginExpansion
{\displaystyle\int\limits_{0}^{2-x-y}}
%EndExpansion
\frac{3}{2}\sqrt{(1-x)(1-y)(1-z)(x+y+z-1)}\,dz\,dy\,dx\\
&  =4\left(  \frac{17\pi}{525}-\frac{\pi^{2}}{160}\right)
\end{align*}
and the mean square area is similarly $4^{2}/560=1/35$.

As far as we know, no one has previously determined the exact density for the
area of a triangle or a cyclic quadrilateral created via broken sticks. \ From
such an expression would come a numerical estimate of the median area
($50\%$-tile), obtained via a single integration. \ We succeed in finding the
density for triangles, but unfortunately not for quadrilaterals. \ Even better
would be an exact cumulative distribution function -- allowing us to avoid the
integration -- and, surprisingly, this too is possible.

\section{PDF\ for Triangle Area}

We work with $z=\operatorname*{area}^{2}$ for now, returning to $\sqrt
{z}=\operatorname*{area}$ at the conclusion. \ The system of equations%
\[%
\begin{array}
[c]{ccc}%
z=(1-x)(1-y)(x+y-1), &  & w=y
\end{array}
\]
has two solutions:%
\[%
\begin{array}
[c]{ccc}%
x=1-\dfrac{w}{2}\pm\dfrac{\sqrt{(1-w)w^{2}-4z}}{2\sqrt{1-w}}, &  & y=w
\end{array}
\]
and the map $(x,y)\mapsto(z,w)$ has absolute Jacobian determinant%
\[
\left\vert
\begin{array}
[c]{cc}%
-(1-y)(2x+y-2) & -(1-x)(x+2y-2)\\
0 & 1
\end{array}
\right\vert =\sqrt{1-w}\sqrt{(1-w)w^{2}-4z}.
\]
Since the joint density for $(x,y)$ is $2$ and the map to $(z,w)$ is
two-to-one, the joint density for $(z,w)$ is \cite{Papoulis}%
\[
\frac{4}{\sqrt{1-w}\sqrt{(1-w)w^{2}-4z}}=\frac{4}{\sqrt{(w-c)(w-a)(b-w)(1-w)}}%
\]
where $c(z)<0<a(z)<w<b(z)<1$ are the three zeroes of the cubic polynomial
$(1-w)w^{2}-4z$:%
\[%
\begin{array}
[c]{l}%
a(z)=\dfrac{1}{3}+\dfrac{1-i\sqrt{3}}{6}\theta(z)^{-1/3}+\dfrac{1+i\sqrt{3}%
}{6}\theta(z)^{1/3},\medskip\\
b(z)=\dfrac{1}{3}+\dfrac{1+i\sqrt{3}}{6}\theta(z)^{-1/3}+\dfrac{1-i\sqrt{3}%
}{6}\theta(z)^{1/3},\medskip\\
c(z)=\dfrac{1}{3}-\dfrac{1}{3}\theta(z)^{-1/3}-\dfrac{1}{3}\theta(z)^{1/3}%
\end{array}
\]
and%
\[
\theta(z)=-1+54z+6\sqrt{3}\sqrt{-z+27z^{2}}.
\]
It follows that the marginal density for $z$ is \cite{GrdRyzk}%
\[
g(z)=%
%TCIMACRO{\dint \limits_{a}^{b}}%
%BeginExpansion
{\displaystyle\int\limits_{a}^{b}}
%EndExpansion
\frac{4\,dw}{\sqrt{1-w}\sqrt{(1-w)w^{2}-4z}}=\frac{8}{\sqrt{(1-a)(b-c)}%
}K\left[  \frac{(b-a)(1-c)}{(1-a)(b-c)}\right]
\]
where $K$ is the complete elliptic integral of the first kind:%
\[
K[m]=%
%TCIMACRO{\dint \limits_{0}^{1}}%
%BeginExpansion
{\displaystyle\int\limits_{0}^{1}}
%EndExpansion
\frac{d\tau}{\sqrt{1-\tau^{2}}\sqrt{1-m\tau^{2}}}%
\]
(consistent with \textsc{Mathematica}). \ The marginal density for $\sqrt{z}$
is therefore%
\[
f(\zeta)=\frac{d}{d\zeta}\mathbb{P}\left\{  \sqrt{z}<\zeta\right\}  =\frac
{d}{d\zeta}\mathbb{P}\left\{  z<\zeta^{2}\right\}  =2\zeta\cdot g\left(
\zeta^{2}\right)
\]
where $0<\zeta<1/\left(  3\sqrt{3}\right)  $; see Figure 3. \
%TCIMACRO{\FRAME{ftbpFU}{3.0277in}{3.0234in}{0pt}{\Qcb{Area density for
%triangles}}{}{areahist.eps}{\special{ language "Scientific Word";
%type "GRAPHIC";  maintain-aspect-ratio TRUE;  display "USEDEF";
%valid_file "F";  width 3.0277in;  height 3.0234in;  depth 0pt;
%original-width 11.6663in;  original-height 11.6516in;  cropleft "0";
%croptop "1";  cropright "1";  cropbottom "0";
%filename 'areahist.eps';file-properties "XNPEU";}} }%
%BeginExpansion
\begin{figure}[ptb]%
\centering
\includegraphics[
height=3.0234in,
width=3.0277in
]%
{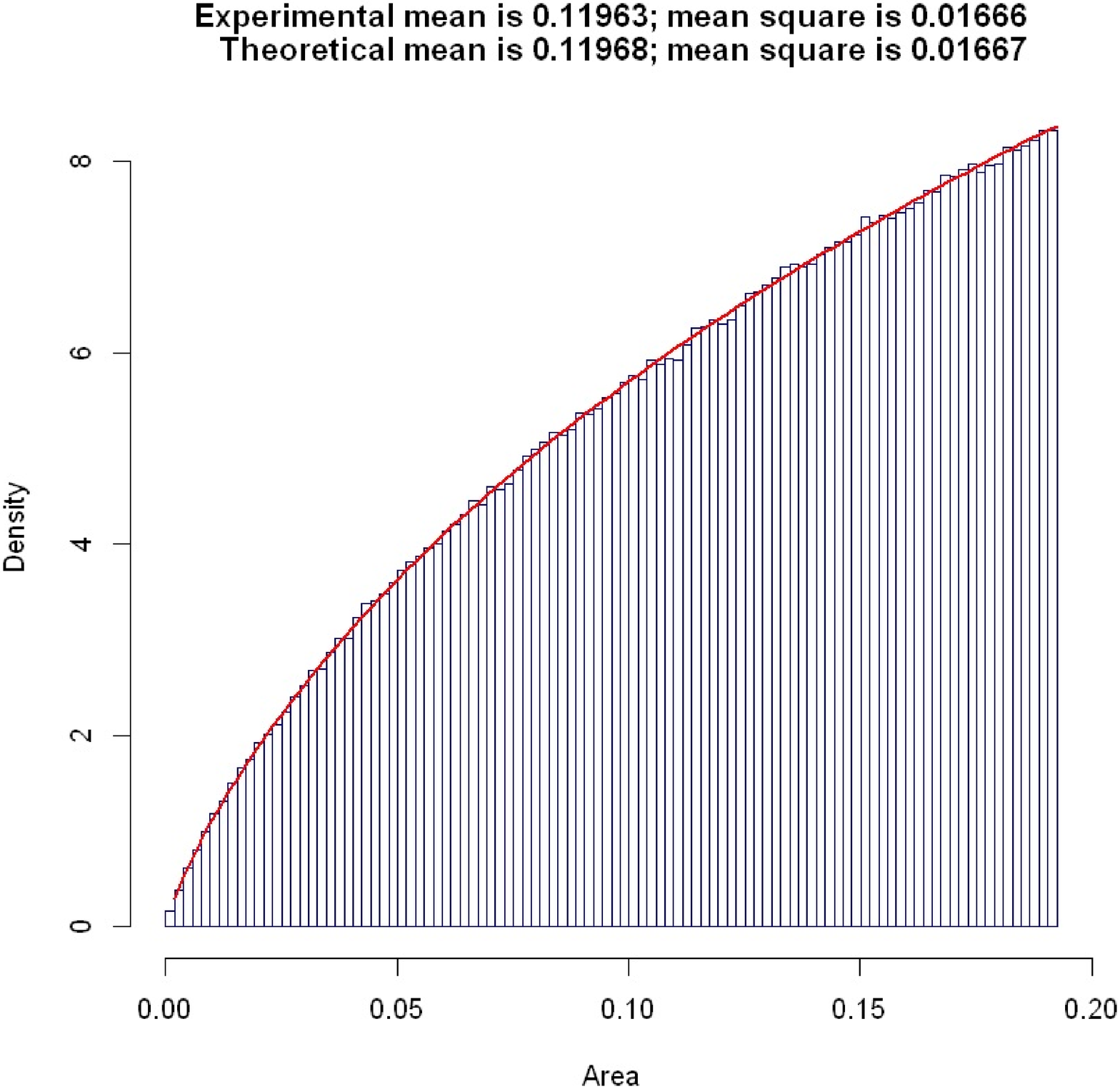}%
\caption{Area density for triangles}%
\end{figure}
%EndExpansion
Numerically solving the equation%
\[%
%TCIMACRO{\dint \limits_{0}^{\mu}}%
%BeginExpansion
{\displaystyle\int\limits_{0}^{\mu}}
%EndExpansion
f(\zeta)\,d\zeta=\frac{1}{2}%
\]
gives the median%
\[
\mu=0.1258338431386510592028005...=(4)(0.0314584607846627648007001...)
\]
to precision limited only by the accuracy of the integration routine. \ A
symbolic antiderivative of $f(\zeta)$ would seem infeasible, at least at first glance.

\section{CDF\ for Triangle Area}

Again, we denote sides by $x$, $y$ and $\operatorname*{area}^{2}$ by $z$.
\ Let $u=2-x-y$ and $v=y-x$, so that $0<u<1$, $-1<v<1$ and%
\[
(u+v)(u-v)(1-u)=(2-2x)(2-2y)\left[  1-(2-x-y)\right]  =4(1-x)(1-y)(x+y-1)=4z
\]
hence%
\[
u^{2}-v^{2}=\frac{4z}{1-u}%
\]
hence%
\[
\left\vert v\right\vert =\sqrt{u^{2}-\frac{4z}{1-u}}=\sqrt{\frac
{(1-u)u^{2}-4z}{1-u}}=q(z,u).
\]
The pair $(u,v)$ is uniform on the domain $\left\{  (u,v):0<u<1\text{ and
}\left\vert v\right\vert <u\right\}  $, a triangle of unit area; thus the
probability that $\operatorname*{area}^{2}$ exceeds $z$ is%
\[%
%TCIMACRO{\dint \limits_{a(z)}^{b(z)}}%
%BeginExpansion
{\displaystyle\int\limits_{a(z)}^{b(z)}}
%EndExpansion
\,%
%TCIMACRO{\dint \limits_{-q(z,u)}^{q(z,u)}}%
%BeginExpansion
{\displaystyle\int\limits_{-q(z,u)}^{q(z,u)}}
%EndExpansion
dv\,du=2%
%TCIMACRO{\dint \limits_{a(z)}^{b(z)}}%
%BeginExpansion
{\displaystyle\int\limits_{a(z)}^{b(z)}}
%EndExpansion
q(z,u)du
\]
where the zeros $c(z)<a(z)<u<b(z)<1$ are exactly as before. \ If $t=\sqrt
{1-u}$, then $u=1-t^{2}$ and $du=-2t\,dt$; it follows that%
\begin{align*}
q(z,u)du  &  =\sqrt{\frac{t^{2}\left(  1-t^{2}\right)  ^{2}-4z}{t^{2}}}\left(
-2t\,dt\right)  =-2\sqrt{\left(  1-t^{2}\right)  ^{2}t^{2}-4z}\,dt\\
&  =-2\sqrt{\left(  1-t^{2}\right)  t-2\sqrt{z}}\,\sqrt{\left(  1-t^{2}%
\right)  t+2\sqrt{z}}\,dt\\
&  =-2\sqrt{\left(  t^{2}-\alpha^{2}\right)  \left(  \beta^{2}-t^{2}\right)
\left(  \gamma^{2}-t^{2}\right)  }\,dt
\end{align*}
where $\beta=\sqrt{1-a}$, $\alpha=\sqrt{1-b}$ and $\gamma=\sqrt{1-c}$. \ From
$1-c>1-a>1-u>1-b>0$, we have $0<\alpha(z)<t<\beta(z)<\gamma(z)$. \ The
preceding argument leading to the formula%
\[
\mathbb{P}\left\{  \operatorname*{area}>\zeta\right\}  =4%
%TCIMACRO{\dint \limits_{\alpha\left(  \zeta^{2}\right)  }^{\beta\left(
%\zeta^{2}\right)  }}%
%BeginExpansion
{\displaystyle\int\limits_{\alpha\left(  \zeta^{2}\right)  }^{\beta\left(
\zeta^{2}\right)  }}
%EndExpansion
\sqrt{\left(  1-t^{2}\right)  ^{2}t^{2}-4\zeta^{2}}\,dt=4J\left(
\zeta\right)
\]
is due to an anonymous student \cite{PACM2}. \ Our only contribution is to
link this with Dieckmann's \cite{Dieckm} integral evaluation:%
\begin{align*}
8J  &  =\beta\sqrt{\gamma^{2}-\alpha^{2}}\left(  \alpha^{2}+\beta^{2}%
+\gamma^{2}\right)  E\left[  \frac{\left(  \beta^{2}-\alpha^{2}\right)
\gamma^{2}}{\left(  \gamma^{2}-\alpha^{2}\right)  \beta^{2}}\right]
+\frac{\alpha^{2}\beta}{\sqrt{\gamma^{2}-\alpha^{2}}}\left(  \alpha^{2}%
+\beta^{2}-5\gamma^{2}\right)  K\left[  \frac{\left(  \beta^{2}-\alpha
^{2}\right)  \gamma^{2}}{\left(  \gamma^{2}-\alpha^{2}\right)  \beta^{2}%
}\right] \\
&  -\frac{\alpha^{2}}{\beta\sqrt{\gamma^{2}-\alpha^{2}}}\left(  \alpha
+\beta-\gamma\right)  \left(  \alpha-\beta-\gamma\right)  \left(  \alpha
-\beta+\gamma\right)  \left(  \alpha+\beta+\gamma\right)  \Pi\left[
\frac{\beta^{2}-\alpha^{2}}{\beta^{2}},\frac{\left(  \beta^{2}-\alpha
^{2}\right)  \gamma^{2}}{\left(  \gamma^{2}-\alpha^{2}\right)  \beta^{2}%
}\right]
\end{align*}
where $E$ and $\Pi$ are complete elliptic integrals of the second and third
kind:%
\[%
\begin{array}
[c]{ccc}%
E[m]=%
%TCIMACRO{\dint \limits_{0}^{1}}%
%BeginExpansion
{\displaystyle\int\limits_{0}^{1}}
%EndExpansion
\dfrac{\sqrt{1-m\tau^{2}}}{\sqrt{1-\tau^{2}}}\,d\tau, &  & \Pi\lbrack n,m]=%
%TCIMACRO{\dint \limits_{0}^{1}}%
%BeginExpansion
{\displaystyle\int\limits_{0}^{1}}
%EndExpansion
\dfrac{d\tau}{\left(  1-n\tau^{2}\right)  \sqrt{1-\tau^{2}}\sqrt{1-m\tau^{2}}%
}.
\end{array}
\]
Solving numerically the equation $8J\left(  \mu\right)  =1$ gives the median
to essentially infinite precision.

\section{Cyclic Quadrilaterals}

On the one hand, arbitrary angles $\alpha$ and $\beta$ in a cyclic
quadrilateral are distributed according to what we call a bivariate tent
density:%
\[
\left\{
\begin{array}
[c]{lll}%
\varphi(\pi-y,x) &  & \text{if }\pi-y<x<y\text{ and }\pi/2<y<\pi\text{,}\\
\varphi(x,y) &  & \text{if }x<y<\pi-x\text{ and }0<x<\pi/2\text{,}\\
\varphi(\pi-x,y) &  & \text{if }\pi-x<y<x\text{ and }\pi/2<x<\pi\text{,}\\
\varphi(y,x) &  & \text{if }y<x<\pi-y\text{ and }0<y<\pi/2
\end{array}
\right.
\]
where%
\[
\varphi(x,y)=\frac{\left[  4\cos(x)-3\cos(2y)-1\right]  \tan\left(
x/2\right)  ^{2}}{2\left[  \sin(x)+\sin(y)\right]  ^{2}\sin(y)^{2}}.
\]
A sketch of the proof is given in Appendix I; a less complicated example
appears in \cite{Fnch1}. \ Clearly $\alpha$ and $\beta$ are uncorrelated yet
dependent. \ The univariate density for $\alpha$ is%
\[
\frac{\psi_{1}(x)-16\psi_{2}(x)\ln(\sin(x)/2)+\psi_{3}(x)\ln(\tan
(x/2))}{16\cos(x)^{3}\sin(x)^{5}}%
\]
where trigonometric polynomials $\psi_{1}$, $\psi_{2}$, $\psi_{3}$ are given
by
\[
\psi_{1}(x)=-25\cos(x)+7\cos(3x)+17\cos(5x)+\cos(7x),
\]%
\[
\psi_{2}(x)=42\cos(x)+19\cos(3x)+3\cos(5x),
\]%
\[
\psi_{3}(x)=378+489\cos(2x)+150\cos(4x)+7\cos(6x)
\]
(Figures 4 and 5).
%TCIMACRO{\FRAME{ftbpFU}{3.8165in}{3.0245in}{0pt}{\Qcb{Bivariate tent densitiy
%on $[0,\pi]\times\lbrack0,\pi]$}}{}{tent1.eps}%
%{\special{ language "Scientific Word";  type "GRAPHIC";
%maintain-aspect-ratio TRUE;  display "USEDEF";  valid_file "F";
%width 3.8165in;  height 3.0245in;  depth 0pt;  original-width 9.263in;
%original-height 7.3258in;  cropleft "0";  croptop "1";  cropright "1";
%cropbottom "0";  filename 'tent1.eps';file-properties "XNPEU";}} }%
%BeginExpansion
\begin{figure}[ptb]%
\centering
\includegraphics[
height=3.0245in,
width=3.8165in
]%
{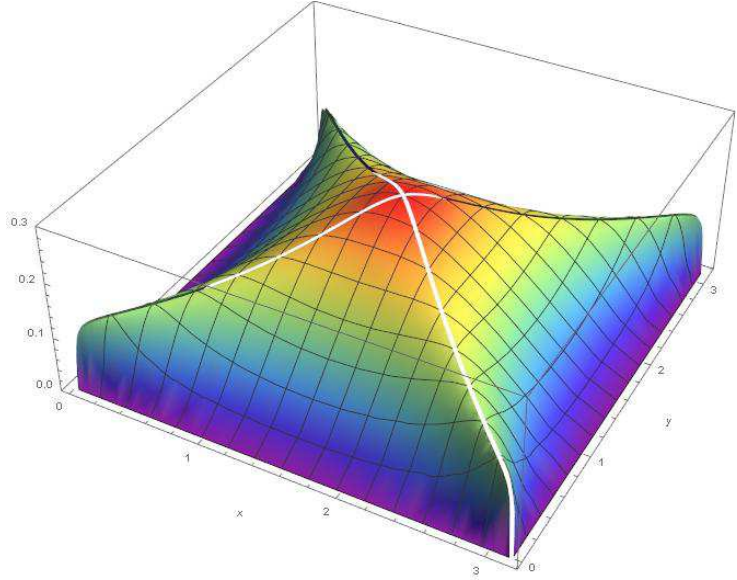}%
\caption{Bivariate tent densitiy on $[0,\pi]\times\lbrack0,\pi]$}%
\end{figure}
%EndExpansion

\
%TCIMACRO{\FRAME{ftbpFU}{3.0277in}{3.0234in}{0pt}{\Qcb{Angle density for cyclic
%quadrilaterals}}{}{qangla.eps}{\special{ language "Scientific Word";
%type "GRAPHIC";  maintain-aspect-ratio TRUE;  display "USEDEF";
%valid_file "F";  width 3.0277in;  height 3.0234in;  depth 0pt;
%original-width 11.6663in;  original-height 11.6516in;  cropleft "0";
%croptop "1";  cropright "1";  cropbottom "0";
%filename 'qangla.eps';file-properties "XNPEU";}} }%
%BeginExpansion
\begin{figure}[ptb]%
\centering
\includegraphics[
height=3.0234in,
width=3.0277in
]%
{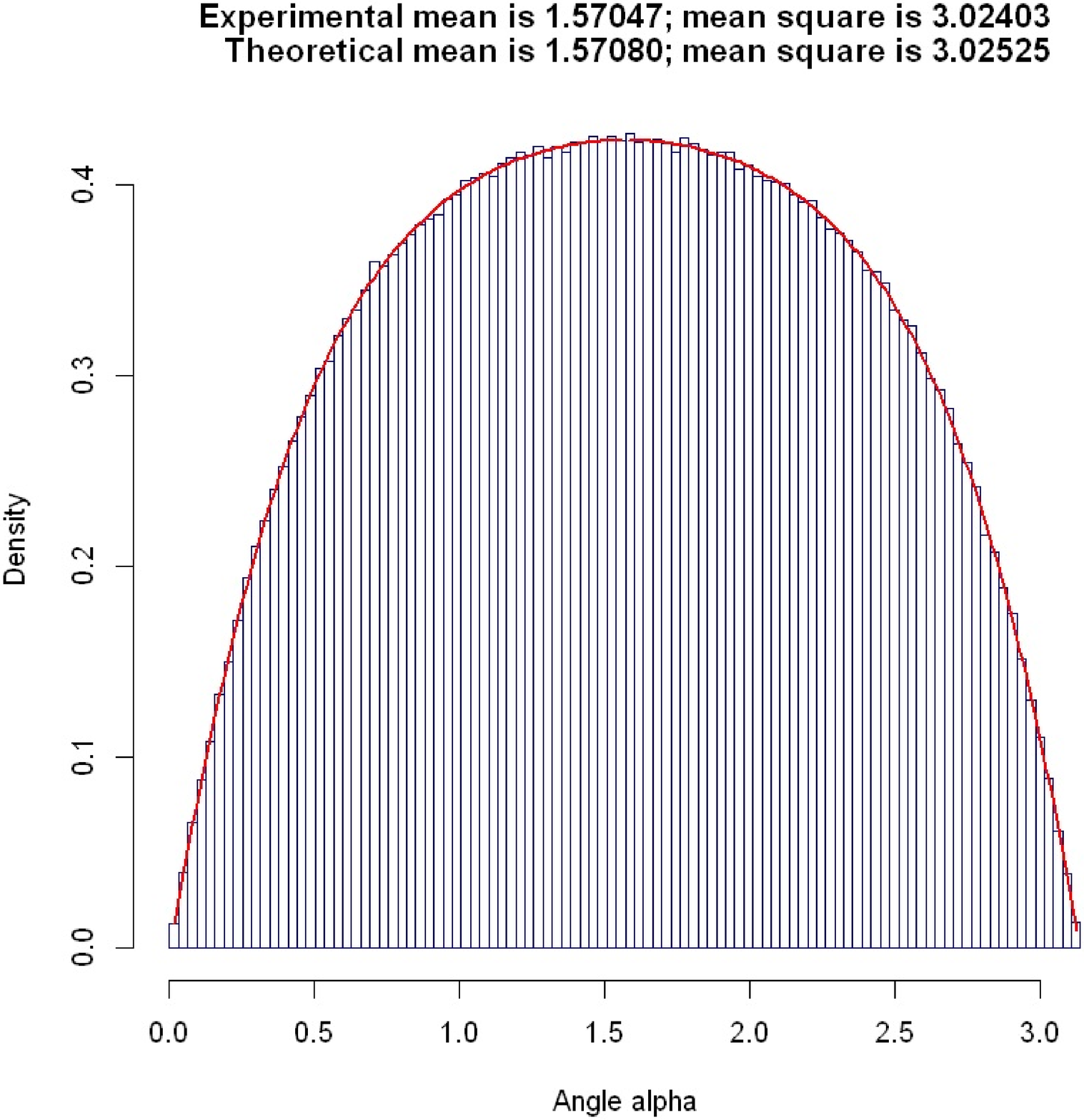}%
\caption{Angle density for cyclic quadrilaterals}%
\end{figure}
%EndExpansion
It follows that $\mathbb{E}(\alpha)=\pi/2$, but a closed-form expression for
\[
\mathbb{E}\left(  \alpha^{2}\right)  =3.0252500344067143300547137...
\]
remains open.

On the other hand, finding the density for area has eluded us -- witness
Appendix II -- and computer simulation suggests that it is approximately
linear (Figure 6). \
%TCIMACRO{\FRAME{ftbpFU}{3.0277in}{3.0234in}{0pt}{\Qcb{Area density for cyclic
%quadrilaterals}}{}{qarea.eps}{\special{ language "Scientific Word";
%type "GRAPHIC";  maintain-aspect-ratio TRUE;  display "USEDEF";
%valid_file "F";  width 3.0277in;  height 3.0234in;  depth 0pt;
%original-width 11.6663in;  original-height 11.6516in;  cropleft "0";
%croptop "1";  cropright "1";  cropbottom "0";
%filename 'qarea.eps';file-properties "XNPEU";}} }%
%BeginExpansion
\begin{figure}[ptb]%
\centering
\includegraphics[
height=3.0234in,
width=3.0277in
]%
{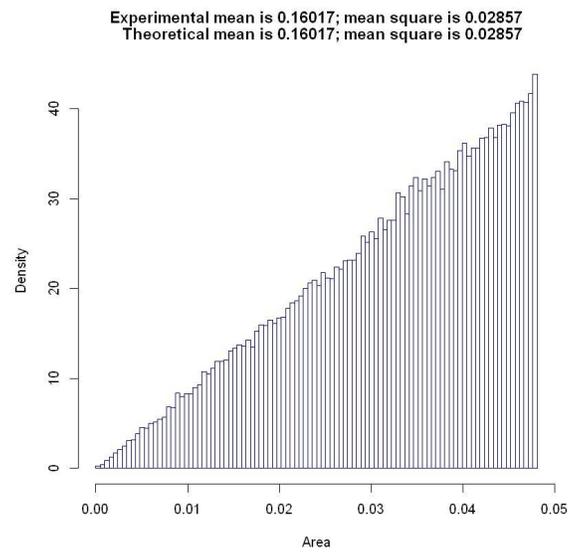}%
\caption{Area density for cyclic quadrilaterals}%
\end{figure}
%EndExpansion
We nonrigorously estimate the median area to be $0.1696...$ via such
experimentation. \ Can this be calculated to high precision?

Of all quadrilaterals with sides $a$, $b$, $c$, $d$ in this order, there is a
unique one with maximal area, the cyclic quadrilateral \cite{Niven, Peter}.
\ The natural analog of this theorem to $n$-gons for $n\geq5$ is true
\cite{Jnkvc, King}.

When breaking a stick in $n-1$ places at random, the $n$ pieces can be
configured as an $n$-gon with probability $1-n/2^{n-1}$ \cite{DAndrG}.

A concise formula for the area of a cyclic pentagon, generalizing those of
Heron and Brahmagupta, apparently does not exist. \ It is known that
$(4\cdot\operatorname*{area})^{2}$ satisfies a $7^{\text{th}}$ degree
polynomial equation with coefficients involving elementary symmetric functions
$\sigma_{k}$ of squares of sides \cite{Rbbns1, Rbbns2, Rbbns3, Rbbns4, Rbbns5,
Rbbns6, Rbbns7}. \ One of two $7^{\text{th}}$ degree polynomials is satisfied
for cyclic hexagons. \ For $n\geq7$, the equations become inconceivably
lengthy, possessing degree $38$ for cyclic heptagons and octagons, and almost
a million terms when expanding with regard to $\sigma_{k}$.

Unless a theoretical breakthrough occurs, broken sticks will never be fully
understood for large $n$. \ A numerical approach is perhaps mandatory. \ We
wonder if even the mean area (let alone the median area) of a cyclic pentagon
is too much for which to ask.

\section{Appendix I}

The bivariate density for two angles of a triangle is easily obtained in
\cite{Fnch2}; the corresponding work for a cyclic quadrilateral is harder.
\ Let adjacent angles $\alpha_{1}$, $\alpha_{2}$ be opposite angles
$\alpha_{3}=\pi-\alpha_{1}$, $\alpha_{4}=\pi-\alpha_{2}$. \ Let sides $s_{1}$,
$s_{2}$ determine $\alpha_{3}$ and sides $s_{3}$, $s_{4}$ determine
$\alpha_{1}$ (see the picture in \cite{Fnch1}). \ By the Law of Cosines,%
\[
s_{1}^{2}+s_{2}^{2}-2s_{1}s_{2}\cos(\alpha_{3})=s_{3}^{2}+s_{4}^{2}%
-2s_{3}s_{4}\cos(\alpha_{1}),
\]%
\[
s_{2}^{2}+s_{3}^{2}-2s_{2}s_{3}\cos(\alpha_{4})=s_{1}^{2}+s_{4}^{2}%
-2s_{1}s_{4}\cos(\alpha_{2})
\]
hence%
\[
s_{1}^{2}+s_{2}^{2}-s_{3}^{2}-s_{4}^{2}=-2(s_{1}s_{2}+s_{3}s_{4})\cos
(\alpha_{1}),
\]%
\[
s_{2}^{2}+s_{3}^{2}-s_{1}^{2}-s_{4}^{2}=-2(s_{2}s_{3}+s_{1}s_{4})\cos
(\alpha_{2})
\]
hence%
\[
\alpha_{1}=\arccos\left(  \frac{s_{3}^{2}+(2-s_{1}-s_{2}-s_{3})^{2}-s_{1}%
^{2}-s_{2}^{2}}{s_{3}(2-s_{1}-s_{2}-s_{3})+s_{1}s_{2}}\right)  ,
\]%
\[
\alpha_{2}=\arccos\left(  \frac{s_{1}^{2}+(2-s_{1}-s_{2}-s_{3})^{2}-s_{2}%
^{2}-s_{3}^{2}}{s_{1}(2-s_{1}-s_{2}-s_{3})+s_{2}s_{3}}\right)
\]
because the stick has length $2$. \ The map $(s_{1},s_{2},s_{3})\mapsto
(\alpha_{1},\alpha_{2},s_{3})$ has absolute Jacobian determinant%
\[
\frac{2(s_{1}+s_{3})(1-s_{3})}{\left[  s_{1}s_{2}+s_{3}(2-s_{1}-s_{2}%
-s_{3})\right]  \left[  s_{2}s_{3}+s_{1}(2-s_{1}-s_{2}-s_{3})\right]  }.
\]
We first rewrite this in terms of $\alpha_{1}$, $\alpha_{2}$, $s_{3}$,
remembering not only $s_{1}>0$, $s_{2}>0$ but also $s_{1}+s_{2}+s_{3}<2$. \ To
do this, perform the substitutions%
\[
s_{1}=1-(1-s_{3})\frac{\tan\left(  \frac{\alpha_{2}}{2}\right)  }{\tan\left(
\frac{\alpha_{1}}{2}\right)  },
\]
\
\[
s_{2}=\frac{-s_{3}\sin\left(  \alpha_{1}+\frac{\alpha_{2}}{2}\right)
+(2-s_{3})\sin\left(  \frac{\alpha_{2}}{2}\right)  }{\sin(\alpha_{1}%
-\frac{\alpha_{2}}{2})+\sin\left(  \frac{\alpha_{2}}{2}\right)  }.
\]
The reciprocal of the determinant is then integrated with respect to $s_{3}$,
with lower limit%
\[
\max\left\{  1-\frac{\tan\left(  \frac{\alpha_{1}}{2}\right)  }{\tan\left(
\frac{\alpha_{2}}{2}\right)  },0\right\}
\]
(from $s_{1}=0$) and upper limit%
\[
\min\left\{  \frac{\sin\left(  \frac{\alpha_{2}}{2}\right)  }{\cos\left(
\frac{\alpha_{1}}{2}\right)  \sin\left(  \frac{\alpha_{1}+\alpha_{2}}%
{2}\right)  },\frac{\cos\left(  \frac{\alpha_{1}}{2}\right)  }{\sin\left(
\frac{\alpha_{2}}{2}\right)  \sin\left(  \frac{\alpha_{1}+\alpha_{2}}%
{2}\right)  }\right\}
\]
(from $s_{2}=0$ and $s_{1}+s_{2}+s_{3}=2$). \ As an example, $s_{2}=0$ when%
\[
s_{3}=\frac{2\sin\left(  \frac{\alpha_{2}}{2}\right)  }{\sin\left(  \alpha
_{1}+\frac{\alpha_{2}}{2}\right)  +\sin\left(  \frac{\alpha_{2}}{2}\right)
}=\frac{\sin\left(  \frac{\alpha_{2}}{2}\right)  }{\cos\left(  \frac
{\alpha_{1}}{2}\right)  \sin\left(  \frac{\alpha_{1}+\alpha_{2}}{2}\right)  }%
\]
since, by the sum-to-product identity,%
\[
\sin\left(  \alpha_{1}+\tfrac{\alpha_{2}}{2}\right)  +\sin\left(
\tfrac{\alpha_{2}}{2}\right)  =2\sin\left(  \frac{\alpha_{1}+\frac{\alpha_{2}%
}{2}+\frac{\alpha_{2}}{2}}{2}\right)  \cos\left(  \frac{\alpha_{1}%
+\frac{\alpha_{2}}{2}-\frac{\alpha_{2}}{2}}{2}\right)  =2\sin\left(
\tfrac{\alpha_{1}+\alpha_{2}}{2}\right)  \cos\left(  \tfrac{\alpha_{1}}%
{2}\right)  .
\]
Likewise, at the end of the derivation,
\[
\frac{\left[  4\cos(\alpha_{1})-3\cos(2\alpha_{2})-1\right]  \tan\left(
\frac{\alpha_{1}}{2}\right)  ^{2}}{32\cos\left(  \frac{\alpha_{1}-\alpha_{2}%
}{2}\right)  ^{2}\sin\left(  \frac{\alpha_{1}+\alpha_{2}}{2}\right)  ^{2}%
\cos\left(  \frac{\alpha_{2}}{2}\right)  ^{2}\sin\left(  \frac{\alpha_{2}}%
{2}\right)  ^{2}}=\frac{\left[  4\cos(\alpha_{1})-3\cos(2\alpha_{2})-1\right]
\tan\left(  \frac{\alpha_{1}}{2}\right)  ^{2}}{2\left[  \sin(\alpha_{1}%
)+\sin(\alpha_{2})\right]  ^{2}\sin(\alpha_{2})^{2}}%
\]
by a double angle formula and since, by the product-to-sum identity,%
\[
2\cos\left(  \tfrac{\alpha_{1}-\alpha_{2}}{2}\right)  \sin\left(
\tfrac{\alpha_{1}+\alpha_{2}}{2}\right)  =\sin\left(  \tfrac{\alpha_{1}%
-\alpha_{2}}{2}+\tfrac{\alpha_{1}+\alpha_{2}}{2}\right)  -\sin\left(
\tfrac{\alpha_{1}-\alpha_{2}}{2}-\tfrac{\alpha_{1}+\alpha_{2}}{2}\right)
=\sin(\alpha_{1})+\sin(\alpha_{2}).
\]
Denote the integral by $I(\alpha_{1},\alpha_{2})$. \ The tent-like appearance
of the surface plot of $I$ suggests necessary simplifications leading to our
formula for the joint density. \ Finally, the details of further integrating
out $\alpha_{2}$ are elaborate and thus omitted.

\section{Appendix II}

We work with $r_{1}=\operatorname*{area}^{2}$. \ The system of equations%
\[%
\begin{array}
[c]{ccccc}%
r_{1}=(1-s_{1})(1-s_{2})(1-s_{3})(s_{1}+s_{2}+s_{3}-1), &  & r_{2}=s_{2}, &  &
r_{3}=s_{3}%
\end{array}
\]
has two solutions:%
\[%
\begin{array}
[c]{ccccc}%
s_{1}=1-\dfrac{r_{2}+r_{3}}{2}\pm\dfrac{\sqrt{(1-r_{2})(1-r_{3})(r_{2}%
+r_{3})^{2}-4r_{1}}}{2\sqrt{1-r_{2}}\sqrt{1-r_{3}}}, &  & s_{2}=r_{2}, &  &
s_{3}=r_{3}%
\end{array}
\]
and the map $(s_{1},s_{2},s_{3})\mapsto(r_{1},r_{2},r_{3})$ has absolute
Jacobian determinant%
\[
\sqrt{1-r_{2}}\sqrt{1-r_{3}}\sqrt{(1-r_{2})(1-r_{3})(r_{2}+r_{3})^{2}-4r_{1}%
}.
\]
Since the joint density for $(s_{1},s_{2},s_{3})$ is $3/2$ and the map to
$(r_{1},r_{2},r_{3})$ is two-to-one, the joint density for $(r_{1},r_{2}%
,r_{3})$ is \cite{Papoulis}%
\[
\frac{3}{\sqrt{1-r_{2}}\sqrt{1-r_{3}}\sqrt{(1-r_{2})(1-r_{3})(r_{2}+r_{3}%
)^{2}-4r_{1}}}=\frac{3}{\sqrt{(r_{3}-c)(r_{3}-a)(b-r_{3})(1-r_{3})(1-r_{2})}}%
\]
where $c(r_{1},r_{2})<a(r_{1},r_{2})<r_{3}<b(r_{1},r_{2})<1$ are the three
zeroes of the cubic polynomial $(1-r_{2})(1-r_{3})(r_{2}+r_{3})^{2}-4r_{1}$
(regarded as a function of $r_{3}$ only). \ What troubles us is that, given
sufficiently small $r_{1}>0$, there is a nonempty interval $\Omega
\subseteq\lbrack0,1]$ for which $r_{2}\in\Omega$ implies $a(r_{1},r_{2})<0$.
\ (As an example, if $r_{1}=0.03$, then $\Omega=[0.4807...,0.8227...]$.)
\ This implies that an integral with respect to $r_{3}$ must possess lower
limit $\max\left\{  0,a(r_{1},r_{2})\right\}  $. \ While this should not
present an obstacle numerically, it does create havoc symbolically. \ To find
exactly the endpoints of $\Omega$, that is, to solve the equation
$a(r_{1},r_{2})=0$ for two values $0<r_{2}^{\prime}<r_{2}^{\prime\prime}<1$
via computer algebra, introduces a complexity roadblock in our stochastic
analysis. \ Such difficulties did not arise in Section 1 because $a(z)$ was
always positive. \ Our hope is that someone else will see a workaround.

\section{Acknowledgements}

I am grateful to Andreas Dieckmann \cite{Dieckm}, who promptly evaluated the
integral containing three quadratic factors at my request. \ It is impressive
to see the recent work of students in \cite{PACM1, PACM2, UIUC1, UIUC2} on
variations of the broken stick; I\ appreciate efforts of their teachers in
keeping the flame of original research alive.

\footnotetext{Copyright \copyright \ 2018 by Steven R. Finch. All rights
reserved.}


\newpage

\begin{thebibliography}{99}                                                                                               %


\bibitem {Dobbs}D. E. Dobbs, The average area of a triangle,
\textit{Mathematics and Computer Education} 21 (1987) 178--181.

\bibitem {Johnson}R. A. Johnson, \textit{Advanced Euclidean Geometry}:
\textit{An Elementary Treatise on the Geometry of the Triangle and the
Circle}, Dover Publications, 1960, pp. 81--85; MR0120538 (22 \#11289).

\bibitem {CxtrGrtz}H. S. M. Coxeter and S. L. Greitzer, \textit{Geometry
Revisited}, Math. Assoc. of Amer., 1967, pp. 56--60; MR3155265.

\bibitem {PACM1}P. A. CrowdMath, The Broken Stick Problem: (ii), MIT
PRIMES/AoPS, 2017, http://artofproblemsolving.com/polymath/mitprimes2017b/p.

\bibitem {Papoulis}A. Papoulis, \textit{Probability, Random Variables, and
Stochastic Processes}, McGraw-Hill, 1965, pp. 187--206; MR0176501 (31 \#773).

\bibitem {GrdRyzk}I. S. Gradshteyn and I. M. Ryzhik, \textit{Table of
Integrals, Series, and Products}, 7$^{\text{th}}$ ed., Elsevier/Academic
Press, 2007, p. 275, s. 3.147; MR2360010 (2008g:00005).

\bibitem {PACM2}P. A. CrowdMath, The Broken Stick Problem: (iii), MIT
PRIMES/AoPS, 2017, http://artofproblemsolving.com/polymath/mitprimes2017b/p.

\bibitem {Dieckm}A. Dieckmann, Table of Indefinite Integrals, Universit\"{a}t
Bonn, http://www-elsa.physik.uni-bonn.de/\symbol{126}dieckman/IntegralsIndefinite/IndefInt.html.

\bibitem {Fnch1}S. R. Finch, Random cyclic quadrilaterals, arXiv:1610.00510.

\bibitem {Niven}I. Niven, \textit{Maxima and Minima Without Calculus}, Math.
Assoc. Amer., 1981, pp. 47--55, 253--256; MR0654149 (83i:52011).

\bibitem {Peter}T. Peter, Maximizing the area of a quadrilateral,
\textit{College Math. J.} 34 (2003) 315-316.

\bibitem {Jnkvc}V. Jankovi\'{c}, On the impossibility of one ruler-and-compass
construction, \textit{Mat. Vesnik} 48 (1996) 73--75; MR1454529 (98c:51025).

\bibitem {King}D. King, Maximum Polygon Area, http://www.drking.org.uk/hexagons/misc/polymax.html.

\bibitem {DAndrG}C. D'Andrea and E. G\'{o}mez, The broken spaghetti noodle,
\textit{Amer. Math. Monthly} 113 (2006) 555--557; MR2231141.

\bibitem {Rbbns1}D. P. Robbins, Areas of polygons inscribed in a circle,
\textit{Discrete Comput. Geom.} 12 (1994) 223--236; MR1283889 (95g:51027).

\bibitem {Rbbns2}D. P. Robbins, Areas of polygons inscribed in a circle,
\textit{Amer. Math. Monthly} 102 (1995) 523--530; MR1336638 (96k:51024).

\bibitem {Rbbns3}F. Miller Maley, D. P. Robbins and J. Roskies, On the areas
of cyclic and semicyclic polygons, Adv. in Appl. Math. 34 (2005) 669--689;
MR2128992 (2006b:51016).

\bibitem {Rbbns4}I. Pak, The area of cyclic polygons: recent progress on
Robbins' conjectures, \textit{Adv. in Appl. Math.} 34 (2005) 690--696;
MR2128993 (2006b:51017).

\bibitem {Rbbns5}D. Svrtan, D. Veljan and V. Volenec, Geometry of pentagons:
from Gauss to Robbins, arXiv:math/0403503.

\bibitem {Rbbns6}P. Pech, Computations of the area and radius of cyclic
polygons given by the lengths of sides, \textit{Automated Deduction in
Geometry}, Lect. Notes in Comput. Sci. 3763, Springer-Verlag, 2006, pp.
44--58; MR2259087 (2008g:51021).

\bibitem {Rbbns7}P. Pech, Computation with pentagons, \textit{J. Geom. Graph.}
12 (2008) 151--160; MR2519392 (2010f:51022).

\bibitem {Fnch2}S. R. Finch, Uniform triangles with equality constraints, arXiv:1411.5216.

\bibitem {UIUC1}L. Kong, L. Lkhamsuren, A. Turner, A. Uppal and A. J.
Hildebrand, Random Points, Broken Sticks, and Triangles, UIUC, 2013,
https://faculty.math.illinois.edu/\symbol{126}hildebr/ugresearch/brokenstick-spring2013report.pdf.

\bibitem {UIUC2}A. Page, Y. Semibratova, Y. Xuan, E. R. Zhang, M. T. Phaovibul
and A. J. Hildebrand, The Broken Stick Problem in Higher Dimensions, UIUC,
2015, https://faculty.math.illinois.edu/\symbol{126}hildebr/ugresearch/Hildebrand-Calculus-Spring2015-report.pdf.

%

\begin{tabular}
[c]{lll}
& Steven Finch & \\
& MIT Sloan School of Management & \\
& Cambridge, MA, USA & \\
& \textit{steven\_finch@harvard.edu} &
\end{tabular}

\end{thebibliography}
\end{document}